\documentclass[12pt,  a4paper]{article}
\usepackage{amssymb}
\usepackage{amsfonts, amsmath, amssymb, amsthm}
\usepackage{latexsym, amscd}
\usepackage{amsbsy}
\usepackage{easymat}

\newcommand{\bcen}{\begin{center}}     \newcommand{\ecen}{\end{center}}
\newcommand{\bay}{\begin{array}}      \newcommand{\eay}{\end{array}}
\newcommand{\beq}{\begin{eqnarray*}}      \newcommand{\eeq}{\end{eqnarray*}}

\def\d{\underline{d}}

\def\dim{\mbox{dim}}

\def\\{\\[1.2ex]}

\begin{document}

\title{\large {\bf Belitskii's canonical forms of linear dynamical systems \footnote{Supported by NSFC (No.11371186, 11571341 and 11301161) and dedicated to Professor Yingbo Zhang on the occasion 70th birthday.}}}

\author{\large Yuan Chen, Liujie Nie and Yunge Xu\footnote{Corresponding author: xuy@hubu.edu.cn (Y.G. Xu)}}

\date{\footnotesize Faculty of Mathematics and Statistics, Hubei Key Laboratory of Applied Mathematics, Hubei
University, Wuhan 430062, P.R. China}

\maketitle

\begin{center}

\begin{minipage}{12cm}

{\bf Abstract}: In the note,  all indecomposable canonical forms of linear systems with dimension less than or equal to $4$ are determined based on Belitskii's algorithm. As an application, an effective way to calculate dimensions of equivalence classes of linear systems is given by using Belitskii's canonical forms.

\medskip


\end{minipage}

\end{center}

\section{Introduction}
Let $\mathbb{R}$ and $\mathbb{C}$ be the fields of real and complex numbers respectively. A linear time-invariant dynamical system  is
denoted by $
\begin{MAT}{|c|c}
\first-
A&B\\4
C&
\addpath{(2,1,4)u}
\addpath{(0,0,4)r}\\
\end{MAT}
$ whose dynamical interpretation is given by
a system of ordinary differential equations
$x'(t) =Ax(t)+Bu(t)$ and $y(t)=Cx(t)$
with $t\in \mathbb{R}$ and $x(t)\in\mathbb{C}^n, u(t)\in\mathbb{C}^m$ and  $y(t)\in\mathbb{C}^l$ respectively \cite{Ta}. Such a system is associated with a representation $V$ of the so-called system quiver $Q$
$$~~~~~~~~~~~~~~~~~~~~~~~\setlength{\unitlength}{1mm}
\begin {picture}(35,10)
\put(-15,1){$\bullet$}  \put(0,1){$\bullet$} \put(15,1){$\bullet$}
\put(-9,3){$\alpha$}  \put(0,9){$\beta$} \put(7,3){$\gamma$}
\put(-15,-3){1} \put(0,-3){2} \put(15,-3){3} \put(-13,2){\vector
(1,0){11}} \put(2,2){\vector (1,0){11}} \put(1,6){\circle{5}}
\put(1.5,3.5){\vector (-1,0){1}}
\end {picture}
\mbox{ with }\begin{array}{l}V_1=\mathbb{C}^m, V_2=\mathbb{C}^n, V_3=\mathbb{C}^l\\
V_{\alpha}=B, V_{\beta}=A, V_{\gamma}=C.\end{array}
$$
This connection was originally noticed by Hazewinkel in \cite{Ha}, and the regular points in system spaces are characterized in \cite{HL} by using representation theory of quivers (see \cite{ARS} for representation theory of quivers).

 Recall that two systems $\begin{MAT}{|c|c}
\first- A&B\\4 C&
\addpath{(2,1,4)u}\addpath{(0,0,4)r}\\\end{MAT}$ and $\begin{MAT}{|c|c}
\first- A'&B'\\4 C'&
\addpath{(2,1,4)u}\addpath{(0,0,4)r}\\\end{MAT}$ are {\it equivalent} if there exist nonsingular matrices $X\in \mathbb{C}^{m\times m},Y\in \mathbb{C}^{n\times n},Z\in \mathbb{C}^{l\times l}$ such that $A'=YAY^{-1}, B'=YBX^{-1}$ and $C'=ZCY^{-1}$, which can also be reformulated in terms of block matrices as follows:
\begin{equation}\left[\begin{MAT}{cc}
A'&B'\\
C'& 0 \\ \end{MAT}\right]=
\left[\begin{MAT}{cc}
Y& \\ &Z \\ \end{MAT}\right]
\left[\begin{MAT}{cc}
A&B\\
C&0
\\ \end{MAT}\right]
\left[\begin{MAT}{cc}
Y& \\ &X
\\ \end{MAT}\right]^{-1}.\end{equation}
Various canonical forms of linear systems have been made in the literature \cite{BG,CVM, He,HP2,Po,PW}.
In this note we will use the Belitskii's algorithm to reduce each system to its canonical form  such that two systems yield the same canonical form if and only if they are equivalent. Roughly speaking, we first reduce the block $A$ to its Jordan normal form $A'=J_A$, and then reduce the blocks $B$ and $C$ in turn using the stabilizer groups of $A'$, and $B',\,\,A'$ respectively. At each stage, the blocks of matrices are repartitioned and one  restricts the transformation group to the stabilizer group of the already reduced blocks. The process ends in a finite number of steps, and produces a canonical form; see \cite{Be}, \cite[Sec. 2.3]{Se1} or \cite{CXL}  for details.

Note that the system quiver $Q$ is a wild quiver \cite{HP1} and thus it seems impossible to give a complete classification of representations of $Q$ and hence linear dynamical systems. However, for some small dimensions, we can give such a complete classification by determining all indecomposable canonical forms of systems $\begin{MAT}{|c|c}
\first- A&B\\4 C&
\addpath{(2,1,4)u}\addpath{(0,0,4)r}\\\end{MAT}$ for dimension $|\d|=m+n+l\leqslant 4$ under the indicated similarity (1) via Belitskii's algorithm.
When $|\d|\geqslant 5$, there are  much more indecomposable canonical forms and thus it is too difficult to give a complete classification of equivalence classes of linear systems.

The rest of the paper is organized as follows.  Section 2  is devoted to calculating all indecomposable canonical forms of linear dynamical systems of dimension $m+n+l\leqslant 4$. In Section 3, the dimension of  equivalence classes of the linear dynamical system is determined via Belitskii's algorithm.

\section{Belitskii's canonical forms of dimension $m+n+1\leqslant 4$}
We will apply Belitskii's algorithm
to compute all indecomposable canonical forms for $|\underline{d}|=m+n+l\leqslant 4$ in this section.

Note that when $m=l=0$, the system $\begin{MAT}{|c|c}
\first- A&B\\4 C&
\addpath{(2,1,4)u}\addpath{(0,0,4)r}\\\end{MAT}$ degenerates to a single matrix $A$ and its Belitskii's canonical form is just the Jordan normal form. So we may assume that $n\ne 0$, and $m$ or $l$ is nonzero, and thus $|\d|\geqslant 2$.

Recall that a system is indecomposable if it can not be written as a direct sum of two non-trivial systems. It is well known that  a system is indecomposable if and only if its endomorphism ring is local. A numerical criterion of indecomposability   was given  in terms of links \cite{XZ}. The system  $
\begin{MAT}{|c|c}
\first-
A&B\\4
C&
\addpath{(2,1,4)u}
\addpath{(0,0,4)r}\\
\end{MAT}
$ of size vector $\d=(m,n,l)$ is indecomposable if and only if the number of links appearing in its canonical form is $|\d|-1=m+n+l-1$.
Since each system is equivalent to a direct sum of finitely many indecomposable systems, it suffices to list all Belitskii's canonical forms of indecomposable systems for each dimension vector $\d=(m,n,l)$.

\textbf{Theorem 1} For $|\d|\leqslant 4$,  the  indecomposable canonical forms of linear systems are as follows
$$
\begin{array}{|c|c|}
\hline
\d & \mbox{indecomposable canonical forms of linear systems of  $|\d|$}\\ \hline
2& \begin{MAT}{|c|c|} \first- \lambda &1\\4 \end{MAT}\;,\;
\begin{MAT}{|c|} \first- \lambda \\4 1\\4 \end{MAT}  \\ \hline
3&
\begin{MAT}{|c|c}
\first- \lambda&1\\4 1&
\addpath{(2,1,4)u}\addpath{(0,0,4)r}\\\end{MAT}\;,\;

\begin{MAT}{|cc|c|}\first- \lambda & 1 &\emptyset \\ &\lambda &1 \\4 \end{MAT}\;,\;

\begin{MAT}{|cc|c|}\first- \lambda & 1 &1 \\ &\lambda &0 \\4 \end{MAT}\;,\;

\begin{MAT}{|cc|c|}\first- \lambda & &1 \\ &\mu &1 \\4 \end{MAT}\;,\;

\begin{MAT}{|cc|}\first-\lambda &1\\ &\lambda\\4 1&\emptyset\\4 \end{MAT}\;,\;

\begin{MAT}{|cc|}\first-\lambda &1\\ &\lambda\\4 0&1\\4 \end{MAT}\;,\;

\begin{MAT}{|cc|}\first-\lambda &\\ &\mu\\4 1&1\\4 \end{MAT}
\\[1ex] \hline
&
\begin{MAT}{|cc|c}
\first- \lambda &1&\emptyset\\ &\lambda&1\\4  1 & \mu &
\addpath{(3,2,4)u}\addpath{(3,1,4)u}\addpath{(0,0,4)r}\addpath{(1,0,4)r} \\\end{MAT}\;,\;

\begin{MAT}{|cc|c}
\first- \lambda &1&\emptyset\\ &\lambda&1\\4  0 & 1 &
\addpath{(3,2,4)u}\addpath{(3,1,4)u}\addpath{(0,0,4)r}\addpath{(1,0,4)r} \\ \end{MAT}\;,\;

\begin{MAT}{|cc|c}
\first- \lambda &1&1\\ &\lambda&0\\4  1 & \emptyset &
\addpath{(3,2,4)u}\addpath{(3,1,4)u}\addpath{(0,0,4)r}\addpath{(1,0,4)r} \\ \end{MAT}\;,\;

\begin{MAT}{|cc|c}
\first- \lambda &1&1\\ &\lambda&0\\4  0 & 1 &
\addpath{(3,2,4)u}\addpath{(3,1,4)u}\addpath{(0,0,4)r}\addpath{(1,0,4)r} \\ \end{MAT}\;,\;

\begin{MAT}{|cc|c}
\first- \lambda &&1\\ &\mu&1\\4  1 & \nu &
\addpath{(3,2,4)u}\addpath{(3,1,4)u}\addpath{(0,0,4)r}\addpath{(1,0,4)r} \\ \end{MAT}\;,\;

\begin{MAT}{|cc|c}
\first- \lambda &&1\\ &\mu&1\\4  0 & 1 &
\addpath{(3,2,4)u}\addpath{(3,1,4)u}\addpath{(0,0,4)r}\addpath{(1,0,4)r} \\ \end{MAT}\;,\;

\begin{MAT}{|cc|c}
\first- \lambda &&0\\ &\mu&1\\4  1 & 1 &
\addpath{(3,2,4)u}\addpath{(3,1,4)u}\addpath{(0,0,4)r}\addpath{(1,0,4)r} \\ \end{MAT}\;,\;

\begin{MAT}{|cc|c}
\first- \lambda &&1\\ &\mu&0\\4  1 & 1 &
\addpath{(3,2,4)u}\addpath{(3,1,4)u}\addpath{(0,0,4)r}\addpath{(1,0,4)r} \\ \end{MAT}\;,\;

\\[3ex]
4&
\begin{MAT}{|ccc|c|}\first-\lambda &1&&\emptyset\\ &\lambda&1&\emptyset\\ &&\lambda&1\\4 \end{MAT}\;,\;

\begin{MAT}{|ccc|c|}\first-\lambda &1&&\emptyset\\ &\lambda&1&1\\ &&\lambda&0\\4 \end{MAT}\;,\;

\begin{MAT}{|ccc|c|}\first-\lambda &1&&1\\ &\lambda&1&0\\ &&\lambda&0\\4 \end{MAT}\;,\;

\begin{MAT}{|ccc|c|}\first-\lambda &1&&\emptyset\\ &\lambda&&1\\ &&\mu &1\\4 \end{MAT}\;,\;

\begin{MAT}{|ccc|c|}\first-\lambda &1&&1\\ &\lambda&&0\\ &&\mu &1\\4 \end{MAT}\;,\;

\begin{MAT}{|ccc|c|}\first-\lambda & &&1\\ &\mu&&1\\ &&\nu &1\\4 \end{MAT}\;,\;

 \\[3ex]
&
\begin{MAT}{|ccc|}\first-\lambda &1&\\ &\lambda&1\\ &&\lambda\\4 1&\emptyset&\emptyset\\4    \end{MAT}\;,\;

\begin{MAT}{|ccc|}\first-\lambda &1&\\ &\lambda&1\\ &&\lambda\\4 0&1&\emptyset\\4 \end{MAT}\;,\;

\begin{MAT}{|ccc|}\first-\lambda &1&\\ &\lambda&1\\ &&\lambda\\4 0&0&1 \\4 \end{MAT}\;,\;

\begin{MAT}{|ccc|}\first-\lambda &1&\\ &\lambda&\\ &&\mu\\4 1&\emptyset&1\\4 \end{MAT}\;,\;

\begin{MAT}{|ccc|}\first-\lambda &1&\\ &\lambda&\\ &&\mu\\4 0&1&1\\4 \end{MAT}\;,\;

\begin{MAT}{|ccc|}\first-\lambda &&\\ &\mu&\\ &&\nu\\4 1&1&1\\4 \end{MAT}\;,\;

\begin{MAT}{|cc|}\first-\lambda &1\\ &\lambda\\4 1&\emptyset\\ 0&1\\4 \end{MAT}\;,\;

\begin{MAT}{|cc|cc|}
\first-\lambda &1& 1&\emptyset\\ &\lambda& 0& 1\\4  \end{MAT}
\\[1ex]\hline
\end{array}
$$
where the parametres $\lambda$, $\mu$ and $\nu$ take pairwise different values in $\mathbb{C}$.

\textbf{Proof.} Based on the formula (1), we reduce in turn the matrices $A,B,C$ to their simplest forms by Belitskii's algorithm. We consider the dimension vector $|\d|$ case by case.

(1) Case $|\d|=2$. There are only two possibilities: $\d=(1,1,0)$ or  $(0,1,1)$. If $\d=(1,1,0)$, then the indecomposable canonical form is
$\begin{MAT}{|c|c|} \first- \lambda &1\\4 \end{MAT},$ where $\lambda\in\mathbb{C}$.
If $\d=(0,1,1)$, then the indecomposable canonical form is
$\begin{MAT}{|c|} \first- \lambda \\4 1\\4 \end{MAT}$.

(2) Case $|\d|=3$. There are five possibilities: $\d=(1,1,1)$, $(1,2,0)$, $(0,2,1)$,  $(0,1,2)$, or  $(2,1,0)$.

When $\d=(1,1,1)$, the formula (1) can be written as $
\left[\begin{smallmatrix}
y& \\ &z
\end{smallmatrix}\right]
\left[\begin{smallmatrix}
a&b\\
c&0
\end{smallmatrix}\right]=\left[\begin{smallmatrix}
a'&b'\\
c'& 0\end{smallmatrix}\right]
\left[\begin{smallmatrix}
y& \\ &x
\end{smallmatrix}\right]$, where $0\ne x,y,z\in\mathbb{C}$.
This implies that $ya=a'y$, which leads to $a'=a$ by loop reduction, and we denote it by $a'=\lambda$. We also have
$yb=b'x$ and by edge reduction we get $b'=1$. Again by edge reduction from the equation $zc=c'y$ we obtain $c'=1$.
So in this case, the indecomposable canonical form is
$\begin{MAT}{|c|c}
\first- \lambda&1\\4 1&
\addpath{(2,1,4)u}\addpath{(0,0,4)r}\\\end{MAT}.$

When $\d=(1,2,0)$, consider the matrix equation
$\left[\begin{smallmatrix}
Y& \\ &x
\end{smallmatrix}\right]
\left[\begin{smallmatrix}
A&B\\
0&0
\end{smallmatrix}\right]=\left[\begin{smallmatrix}
A'&B'\\
0& 0\end{smallmatrix}\right]
\left[\begin{smallmatrix}
Y& \\ &x
\end{smallmatrix}\right]$
where $A,A'\in\mathbb{C}^{2\times 2}, B,B'\in\mathbb{C}^{2\times 1}$, and $Y$ is of size $2\times 2$. The first reduction equation is $YA=A'Y$, then $A'=\left[\begin{smallmatrix}
\lambda&1\\ &\lambda
\end{smallmatrix}\right]$ or $\left[\begin{smallmatrix}
\lambda&\\ &\mu
\end{smallmatrix}\right]$ (possibly, $\lambda=\mu$). We consider reductions for $B$ and $C$ case by case according to $A'$.
\begin{enumerate}
  \item [1)] If $A'=\left[\begin{smallmatrix}
\lambda&1\\ &\lambda
\end{smallmatrix}\right]$, then we reduce the block $B$ using the stabilizer group of $A'$, which is given by  $Y=\left[\begin{smallmatrix}
y&y_1\\ &y
\end{smallmatrix}\right]$. In this case $B$ is partitioned into $B=\left[\begin{smallmatrix}b_2\\ b_1
\end{smallmatrix}\right]$.
 We next reduce the subblock $b_1$ based on the equation $yb_1=b'_1x$ which is obtained from the equation $YB=B'x$. This yields $b'_1=1$ or $0$ by edge reduction. If $b'_1=1$, then the stabilizer group of $A'$ and $b'_1$ is given by $Y=\left[\begin{smallmatrix}y &y_1\\ &y\end{smallmatrix}\right]$ and $y=x$, which is used to reduce the subblock $b_2$ based on the reduction equation $yb_2+y_1=b'_2x$. This leads to $b'_2=0$ (and we denote it by $b'_2=\emptyset$) by regularization.
 If $b'_1=0$, then the reduction equation $yb_2=b'_2x$ leads to $b'_2=1$. Clearly, the canonical forms $\begin{MAT}{|cc|c|}\first- \lambda & 1 &\emptyset \\ &\lambda &1 \\4 \end{MAT}$ and $\begin{MAT}{|cc|c|}\first- \lambda & 1 &1 \\ &\lambda &0 \\4 \end{MAT}$ are indecomposable.
\item [2)] If $A'=\lambda I_2$, then the next reduction equation is $YB=B'x$ which yields $B'=\left[\begin{smallmatrix}1\\ 0\end{smallmatrix}\right]$
or $B'=\left[\begin{smallmatrix}0\\ 0\end{smallmatrix}\right]$.
Thus we get the canonical forms $\begin{MAT}{|cc|c|}\first- \lambda & 0 &1 \\ &\lambda &0 \\4 \end{MAT}$ and $\begin{MAT}{|cc|c|}\first- \lambda & 0 &0 \\ &\lambda &0 \\4 \end{MAT}$ which are decomposable.
\item [3)] If $A'=\left[\begin{smallmatrix}
\lambda&\\ &\mu
\end{smallmatrix}\right]$ with $\lambda\ne \mu$, then the stabilizer group of $A'$ is given by
$Y=\left[\begin{smallmatrix}
y_1&\\ &y_2
\end{smallmatrix}\right]$. Write $B=\left[\begin{smallmatrix}b_2\\ b_1
\end{smallmatrix}\right]$. Then the reduction equations for $b_1$ and $b_2$ are $y_2b_1=b'_1x$ and $y_1b_2=b'_2x$, respectively. If $b_1\ne 0$ and $b_2\ne 0$, then we have $b'_1=b'_2=1$ by edge reduction and the canonical form is
$\begin{MAT}{|cc|c|}\first- \lambda & &1 \\ &\mu &1 \\4 \end{MAT}$  which is indecomposable. If $b_1=0$ or $b_2=0$, then $b'_1=0$ or $b'_2=0$, and thus we obtain the canonical form
$\begin{MAT}{|cc|c|}\first- \lambda & &1 \\ &\mu &0 \\4 \end{MAT}$ or $\begin{MAT}{|cc|c|}\first- \lambda & &0 \\ &\mu &1 \\4 \end{MAT}$ which is decomposable.
\end{enumerate}

When $\d=(0,2,1)$, we can similarly obtain other three canonical forms $\begin{MAT}{|cc|}\first-\lambda &1\\ &\lambda\\4 1&\emptyset\\4 \end{MAT}$,
$\begin{MAT}{|cc|}\first-\lambda &1\\ &\lambda\\4 0&1\\4 \end{MAT}$ and
$\begin{MAT}{|cc|}\first-\lambda &\\ &\mu\\4 1&1\\4 \end{MAT}$. Note that for dimension vectors  $(2,1,0)$ and $(0,1,2)$, their canonical forms are all decomposable by similar discussions. We complete the proof in this case.

(3) Case $|\d|=4$. There are nine possibilities: $\d=(1,1,2)$, $(1,2,1)$,  $(2,1,1)$, $(1,3,0)$, $(3,1,0)$, $(0, 1, 3)$, $(0,3,1)$, $(2,2,0)$ or  $(0,2,2)$.

When $\d=(1,2,1)$, the reduction equation
$\left[\begin{smallmatrix}
Y& \\ &z
\end{smallmatrix}\right]
\left[\begin{smallmatrix}
A&B\\
C&0
\end{smallmatrix}\right]=\left[\begin{smallmatrix}
A'&B'\\
C'& 0\end{smallmatrix}\right]
\left[\begin{smallmatrix}
Y& \\ &x
\end{smallmatrix}\right]$
yields $YA=A'Y$, and thus $A'=\left[\begin{smallmatrix}
\lambda&1\\ &\lambda\end{smallmatrix}\right]$ or $\left[\begin{smallmatrix}
\lambda&\\ &\mu\end{smallmatrix}\right]$ (possibly $\lambda= \mu$) by loop reduction. Then according to $A'$, we have
\begin{enumerate}
  \item [1)] If $A'=\left[\begin{smallmatrix}
\lambda&1\\ &\lambda\end{smallmatrix}\right]$, then the stabilizer group of $A'$ is given by
$Y=\left[\begin{smallmatrix}
y&y_1\\ &y
\end{smallmatrix}\right]$. Write $B=\left[\begin{smallmatrix} b_2\\ b_1\end{smallmatrix}\right]$
and $C=[c_1\ c_2]$. If $b_1\ne 0$, then the reduction equations
\begin{equation}
yb_1=b'_1x, \quad yb_2+y_1b_1=b'_2x
\end{equation}
lead to $b'_1=1$ by edge reduction and $b'_2=\emptyset$ by regularization. Moreover, if $c_1\ne 0$, then the reduction equations
\begin{equation}
zc_1=c'_1y, \quad zc_2=c'_1y_1+c'_2y
\end{equation}
lead to  $c'_1=1$ by edge reduction and $c'_2=\mu$ by loop reduction. So we have an indecomposable canonical form
$\begin{MAT}{|cc|c}
\first- \lambda &1&\emptyset\\ &\lambda&1\\4  1 & \mu &
\addpath{(3,2,4)u}\addpath{(3,1,4)u}\addpath{(0,0,4)r}\addpath{(1,0,4)r} \\\end{MAT}$\;.
If $c_1=0$ but $c_2\ne 0$, then by (3), we have $c'_1=0$ and $c'_2=1$ by edge reductions, and thus we have an indecomposable canonical form
$\begin{MAT}{|cc|c}
\first- \lambda &1&\emptyset\\ &\lambda&1\\4  0 & 1 &
\addpath{(3,2,4)u}\addpath{(3,1,4)u}\addpath{(0,0,4)r}\addpath{(1,0,4)r} \\ \end{MAT}\;.$
If $b_1=0$ but $b_2\ne 0$, then by (2), we have $b'_1=0$ and $b'_2=1$ by edge reductions. Consider the two cases that $c_1\ne 0$, and $c_1=0$ but $c_2\ne 0$, we have indecomposable canonical forms
$\begin{MAT}{|cc|c}
\first- \lambda &1&1\\ &\lambda&0\\4  1 & \emptyset &
\addpath{(3,2,4)u}\addpath{(3,1,4)u}\addpath{(0,0,4)r}\addpath{(1,0,4)r} \\ \end{MAT}$ and
$\begin{MAT}{|cc|c}
\first- \lambda &1&1\\ &\lambda&0\\4  0 & 1 &
\addpath{(3,2,4)u}\addpath{(3,1,4)u}\addpath{(0,0,4)r}\addpath{(1,0,4)r} \\ \end{MAT}$
respectively.
If $b_1=b_2=0$, then there is no indecomposable canonical form.

  \item [2)] If $A'=\left[\begin{smallmatrix}
\lambda&\\ &\mu\end{smallmatrix}\right]$ with $\lambda\ne \mu$, then
the stabilizer group of $A'$ is given by  $Y=\left[\begin{smallmatrix}
y_1&\\ &y_2
\end{smallmatrix}\right]$. Write $b=\left[\begin{smallmatrix} b_2\\ b_1\end{smallmatrix}\right]$
and $C=[c_1\ c_2]$. Then the  reduction equations $y_2b_1=b'_1x$, $y_1b_2=b'_2x$, $zc_1=c'_1y_1$ and $zc_2=c'_2y_2$
yield the indecomposable canonical forms
$\begin{MAT}{|cc|c}
\first- \lambda &&1\\ &\mu&1\\4  1 & \nu &
\addpath{(3,2,4)u}\addpath{(3,1,4)u}\addpath{(0,0,4)r}\addpath{(1,0,4)r} \\ \end{MAT}$\;, $\begin{MAT}{|cc|c}
\first- \lambda &&1\\ &\mu&1\\4  0 & 1 &
\addpath{(3,2,4)u}\addpath{(3,1,4)u}\addpath{(0,0,4)r}\addpath{(1,0,4)r} \\ \end{MAT}$\;, $\begin{MAT}{|cc|c}
\first- \lambda &&0\\ &\mu&1\\4  1 & 1 &
\addpath{(3,2,4)u}\addpath{(3,1,4)u}\addpath{(0,0,4)r}\addpath{(1,0,4)r} \\ \end{MAT}$\;, $\begin{MAT}{|cc|c}
\first- \lambda &&1\\ &\mu&0\\4  1 & 1 &
\addpath{(3,2,4)u}\addpath{(3,1,4)u}\addpath{(0,0,4)r}\addpath{(1,0,4)r} \\ \end{MAT}$\;.

  \item [3)] If $A'=\left[\begin{smallmatrix}
\lambda&\\ &\lambda\end{smallmatrix}\right]$, then the reduction equation $YB=B'x$ gives
$B'=\left[\begin{smallmatrix}1\\0\end{smallmatrix}\right]$ or $B'=0$.
If $B'=\left[\begin{smallmatrix}1\\0\end{smallmatrix}\right]$, then the stabilizer group of $B',\,A'$ is given by  $\left[\begin{smallmatrix}
z&&\\ &Y&\\ &&x
\end{smallmatrix}\right]$ where $Y=\left[\begin{smallmatrix}
x&y_1\\ &y
\end{smallmatrix}\right]$.
The reduction equations for the block $C=[c_1\ c_2]$ are $zc_1=c'_1x$ and $zc_2=c'_1y_1+c'_2y$, which lead to $C'=[1\ \emptyset]$ or $[0\ 1]$ depending on whether $c_1$ is zero or not. In both cases, the canonical forms are decomposable.
If $B'=0$, then the corresponding canonical form is decomposable since $C'$ contains at most two links.
\end{enumerate}

Combining the above three cases, we get eight indecomposable canonical forms when $\d=(1,2,1)$.

By similar discussions, when $\d=(1,3,0)$, we get six indecomposable canonical forms, the 9th to the 14th ones in the case of $|\d|=4$ in the above table in Theorem 1; when $\d=(0,3,1)$, we have another six  indecomposable canonical forms which are the 15th to the 20th ones; when $\d=(2,2,0)$ or $\d=(0,2,2)$, there are only two indecomposable canonical forms which are given as the last two ones in the table for the case of $|\d|=4$; when $\d=(3,1,0)$ or $\d=(0,1,3)$,  there is no indecomposable canonical form. We complete the proof.

\section{Dimensions of equivalence classes of linear dynamical systems}

Belitskii's canonical forms provide not only spectral information, but also an effective way to calculate dimensions of equivalence classes of linear systems.

For any given system $\begin{MAT}{|c|c}
\first- A&B\\4 C&
\addpath{(2,1,4)u}\addpath{(0,0,4)r}\\\end{MAT}$\;, one can reduce it to Belitskii's canonical form by a finite number of steps of reductions. At each stage, one reduces a subblock $M_k\in \mathbb{C}^{m_k\times n_k}$ of the system  to one of the three forms: $\emptyset, \left[\begin{smallmatrix} 0& I_r\\ 0&0\end{smallmatrix}\right]$ or $W$. Here $W$ is a Weyr matrix obtained by some simultaneous permutations of rows and columns of a Jordan matrix  $J=\oplus_{i=1}^s\oplus_{j=1}^{t_i} J_{q_{ij}}(\lambda_i)$ with $q_{i1}\geqslant q_{i2}\geqslant \cdots \geqslant q_{it_i}\geqslant 1$ for $i=1,2,\ldots,s$. If $M_k=\emptyset$, then we define $\sigma_k(M)=m_kn_k$; if $M_k=\left[\begin{smallmatrix} 0 & I_r \\
0 & 0 \end{smallmatrix}\right]$, then we define $\sigma_k(M)=r(m_k+n_k-r)$; if $M_k=W$, then we define
$
\sigma_k(M)=m_k^2-\sum_{i=1}^s\sum_{j=1}^{t_i}(2j-1)q_{ij}.
$
Now we have

\textbf{Theorem 2} Let $\Sigma=\begin{MAT}{|c|c}
\first- A&B\\4 C&
\addpath{(2,1,4)u}\addpath{(0,0,4)r}\\\end{MAT}$ be a linear dynamical system. Then the dimension of equivalence classes of  $\Sigma$ is $\sum_{k=1}^u\sigma_k$.

\textbf{Proof.} Suppose that $M_1, M_2,\ldots,M_u$ are all the reduced subblocks in the canonical form of $\Sigma$, and $G_i$ is the stabilizer group of $M_1, M_2, \ldots, M_i$. The algebraic group $G=Gl_m\times Gl_n\times Gl_l$ acts on the system space $\mathcal{S}=\mathcal{S}_{m,n,l}=\mathbb{C}^{mn}\times \mathbb{C}^{n^2} \times \mathbb{C}^{ln}$
by conjugacy: $(X,Y,Z)\cdot  (A,B, C)
\mapsto (YAY^{-1},YBX^{-1}, ZCY^{-1})$.
It is well known that
$\dim\mathcal{O}_{\Sigma}+\dim \,\mbox{stab}_G(\Sigma)=\dim G$.

By Belitskii's algorithm, it is not difficult to see that $\sigma_k=\dim G_{k-1}-\dim G_k$
by  a classical result \cite[Thm. 2.1]{DE} and \cite[Lemmas 3-5]{XZ}.
By adding both sides of the above equalities for $k=1,2,\ldots,u$, we obtain
$\dim G-\dim G_u=\sum_{k=1}^u\sigma_k.$

Note that $\mbox{stab}_G(\Sigma)=G_u$ and thus
$\dim\mathcal{O}_{\Sigma}=\dim\, G-\dim\,\mbox{stab}_G(\Sigma)=\dim G-\dim G_u=\sum_{k=1}^u\sigma_k$ as desired.

\textbf{Examples.} (1) For $\Sigma=\begin{MAT}{|cc|c}
\first- \lambda &1&\emptyset\\ &\lambda&1\\4  1 & \mu &
\addpath{(3,2,4)u}\addpath{(3,1,4)u}\addpath{(0,0,4)r}\addpath{(1,0,4)r}
\addpath{(1,0,.)u}\addpath{(2,2,.)r}\\\end{MAT}
$ as in the proof of Theorem 1, we have $\sigma_1=2, \sigma_2=\sigma_3=\sigma_4=1, \sigma_5=0$ and thus the dimension of equivalence classes of $\Sigma$ is $5$ by Theorem 2.

(2) For $|\d|=3$, by Theorem 2, the dimensions of equivalence classes of the seven  linear systems listed in the third row of the table in Theorem 1 are $2,4,4,4,4,4,4$, respectively. They are all maximal dimension (i.e., regular points in system space $\mathcal{S}_{\d}$; see \cite{HL}).

\end{document}